\documentclass[amsmath,amssymb,floatfix,preprint,showpacs,superscriptaddress]{revtex4}

\usepackage{graphicx,epsfig}
\usepackage{dcolumn}
\usepackage{bm}
\usepackage{amsmath}
\usepackage{amsthm}

\begin{document}


\title{Semi-Markov Graph Dynamics}

\author{Marco Raberto}
\email{marco.raberto@unige.it}
\homepage{http://ideas.repec.org/e/pra66.html}
\affiliation{Dipartimento di Ingegneria Biofisica ed Elettronica,
Universit\`a degli Studi di Genova, Italy}

\author{Fabio Rapallo}
\email{fabio.rapallo@mfn.unipmn.it}
\homepage{people.unipmn.it/~rapallo} \affiliation{Dipartimento di
Scienze e Tecnologie Avanzate, Universit\`a del Piemonte Orientale
``Amedeo Avogadro'', Viale Michel 11, 15121 Alessandria, Italy}

\author{Enrico Scalas}
\email{enrico.scalas@mfn.unipmn.it}
\homepage{people.unipmn.it/~scalas} \affiliation{Dipartimento di Scienze e Tecnologie Avanzate,
Universit\`a del Piemonte Orientale ``Amedeo Avogadro'', Viale
Michel 11, 15121 Alessandria, Italy and Basque Center for Applied Mathematics, Bizkaia Technology Park, Building 500
48160 Derio, Spain}

\date{\today}



\begin{abstract}

In this paper, we outline a model of graph (or network) dynamics based on two ingredients. The first
ingredient is a Markov chain on the space of possible graphs. The second ingredient is a semi-Markov
counting process of renewal type. The model consists in subordinating the Markov chain to the semi-Markov
counting process. In simple words, this means that the chain transitions occur at random time instants called epochs.
The model is quite rich and its possible connections with algebraic geometry
are briefly discussed. Moreover, for the sake of simplicity, we focus on the space of undirected graphs with
a fixed number of nodes. However, in an example, we present an interbank market model where it
is meaningful to use directed graphs or even weighted graphs.

\end{abstract}

\maketitle

\section{Introduction}
\label{intro}

The publication of {\em Collective dynamics of `small world'
networks} by Watts and Strogatz \cite{watts1998} gave origin to a
plethora of papers on network structure and dynamics. The history
of this scientific fashion is well summarized by Rick Durrett
\cite{durrett}:
\begin{quote}
The theory of random graphs began in the late 1950s in several papers by Erd\"os and R\'enyi. In the late
twentieth century, the notion of six degrees of separation, meaning that any two people on the planet
can be connected by a short chain of people who know each other, inspired Strogatz and Watts
\cite{watts1998} to define
the small world random graph in which each side is connected to $k$ close neighbors, but also
has long-range connections. At about the same time, it was observed in human social and sexual networks
and on the Internet that the number of of neighbors of an individual or computer has a power law
distribution. This inspired Barab\'asi and Albert \cite{barabasi1999} to define the preferential attachment model, which
has this properties. These two papers have led to an explosion of research. While this literature is
extensive, many of the papers are based on simulations and nonrigorous arguments.
\end{quote}
Incidentally, the results of Watts and Strogatz were inspired by
the empirical and theoretical work by Milgram \cite{milgram1967}
and Granovetter \cite{granovetter1973} back in the 1960s and
1970s; similarly, the preferential attachment model by Barab\'asi
and Albert is closely related to the famous 1925 paper by Yule
\cite{yule1925} as well as to a celebrated work by Herbert Simon
published in 1955 \cite{simon1955} (see also chapters 8 and 9 in
reference \cite{garibaldi10} for a recent analysis on Simon's
results). This body of literature is partially reviewed in
Durrett's book \cite{durrett} as well as in a popular science book
written by Barab\'asi \cite{barabasi03}.

It might be interesting to understand why this scientific fashion
was born and how. On this respect, we can quote Wikipedia's
article (as retrieved on 4 May 2011) on Milgram's experiment in
popular culture \cite{wiki1}:
\begin{quote}
Social networks pervade popular culture in the United States and elsewhere. In particular, the notion
of six degrees has become part of the collective consciousness. Social networking websites such as
Friendster, MySpace, XING, Orkut, Cyworld, Bebo, Facebook, and others have greatly increased the
connectivity of the online space through the application of social networking concepts. The
\textquotedblleft Six Degrees'' Facebook application calculates the number of steps between any two
members. [$\ldots$]
\end{quote}
In other words, the social character of human beings combined with the hyper-simplification (trivialization)
of some results promoted by leading science journals might have triggered interest in social
networkology also outside scientific circles.
Moreover, the emergence of social networks in the Internet has indeed made some tools developed by
networkologists profitable. However, a deeper analysis by sociologists and historians of science will be necessary
to falsify or corroborate such hypotheses.

In this paper, we pay our tribute to this fashion, but we slightly
depart from the bulk of literature on social network dynamics.
First of all we consider time evolution also in {\em continuous
time} and not only in discrete time. As the reader will see, this
will be enough to give rise to interesting non stationarities as
well as to non-trivial ergodic behavior. Moreover, to begin with a
simple situation, we will be concerned with {\em undirected
graphs} whose number of nodes $M$ does not change in time.
These restrictions can be easily overcome and, indeed, in the following,
an example with {\em directed graphs} will be presented.
The dynamic variable will be the {\em topology} of the
graph. This approach is motivated by the following considerations.
Social networks are intrinsically volatile. You can be in contact
with someone for a finite time (at a meeting, during a phone call,
etc.), but never meet this person again in the future. This
interaction may or may not have effects on your future actions. If
memory is not a major issue, the new configuration of the graph
will only depend on the previous configuration. Memory {\em is}
indeed an issue, but again, to simplify the analysis, we will
consider a semi-Markov dynamics on the state space of all the
possible graphs with $M$ nodes. It is already quite rich. Incidentally,
notice that, except for the case of infinite memory, finite memory
processes in discrete time are Markov chains.

The dynamics will be defined by a Markov chain subordinated to a
generic counting process. Similar models have been around for many
years. They were (and are) commonly used in engineering and
decision analysis and, on this point, the interested reader can
consult the monograph by Howard \cite{howard1971}.

In this framework, it is often assumed that the waiting times
between consecutive events do follow the exponential distribution,
so that the corresponding counting process is Poisson. Indeed,
many counting processes with non-stationary and non-independent
increments converge to the Poisson process after a transient. If
these counting processes are renewal, i.e. inter-arrival times
$\{J_i\}_{i=1}^\infty$ are independent and identically distributed
(iid) random variables, it is sufficient to assume that the
expected value of these inter-arrival times is finite. However,
recently, it has been shown that heavy-tailed distributed
interarrival times (for which $\mathbb{E}(J_i)= \infty$) play an
important role in human dynamics
\cite{scalas2004,scalas2006,barabasi2010}. After defining the
process in Section II, we will present two introductory examples
in Section III and a detailed model of interbank market in Section
IV.

\section{Theory}

This section begins with the definition of the two basic ingredients of our model, namely
\begin{enumerate}
\item a {\em discrete-time Markov chain} on the finite set of $2^{M(M+1)/2}$ undirected graphs
with $M$ vertices (nodes), and
\item a counting process $N(t)$ for the point process corresponding to a renewal process.
\end{enumerate}
The rest of the section is devoted to the definition of the basic model class.

\subsection{Ingredient 1: a Markov chain on graphs}

Consider an undirected graph $\mathcal{G}_M = (V_M,E)$ where $V_M$
represents a set of $M$ vertices (nodes) and $E$ the corresponding
set of edges. Any such undirected graph can be represented by a
symmetric $M \times M$ adjacency matrix
$\mathbb{A}_{\mathcal{G}_M}$, or simply $\mathbb{A}$, with entries
$A_{i,j} = A_{j,i} = 1$ if vertices $i$ and $j$ are connected by
an edge and $A_{i,j} = A_{j,i} = 0$ otherwise. Note that algebraic
graph theory using linear algebra leads to many interesting
results relating spectral properties of adjacency matrixes to the
properties of the corresponding graphs \cite{agt,biggs93}. For
instance, the matrix
$$
\mathbb{A} = \left( \begin{array}{ccc}
0 & 1 & 1\\
1 & 0 & 1 \\
1 & 1 & 0
\end{array} \right)
$$
corresponds to a graph where there are no self-connections and
each vertex is connected to the other two vertices. As mentioned
above, for a given value of $M$ there are $2^{M(M+1)/2}$ possible
graphs. To see that, it is sufficient to observe that the $M$
diagonal entries can assume either value 1 or value 0 and the same
is true for the $M(M-1)/2$ upper diagonal entries. Now, denote by
$G_M$ the set of $2^{M(M+1)/2}$ undirected graphs with $M$ nodes.
Consider a sequence of random variables $X_1, \ldots, X_n$
assuming values in $G_M$. This becomes our {\em state space}, and
the set of $n$ random variables is a finite stochastic process.
Its full characterization is in term of all finite dimensional
distributions of the following kind (for $1 \leq m \leq n$)
\cite{billingsley}
\begin{equation}
\label{fidi1}
p_{X_1,\ldots,X_m} (x_1, \ldots, x_m) =
\mathbb{P}( X_1 = x_1, \ldots, X_m = x_m),
\end{equation}
where $\mathbb{P}(\cdot)$ denotes the probability of an event
with the values $x_i$ running on all the possible graphs $\mathcal{G}_M$ of $G_M$. The finite dimensional
distributions defined in equation (\ref{fidi1}) obey the two compatibility conditions of Kolmogorov
\cite{billingsley}, namely a symmetry condition
\begin{equation}
\label{sym}
p_{X_1,\ldots,X_m} (x_1, \ldots, x_m) = p_{X_{\pi_1},\ldots,X_{\pi_m}}(x_{\pi_1},\ldots,x_{\pi_m})
\end{equation}
for any permutation $(\pi_1, \ldots, \pi_m)$ of the $m$ random
variables (this is a direct consequence of the symmetry property
for the intersection of events) and a second condition
\begin{equation}
\label{cond2}
p_{X_1,\ldots,X_m} (x_1, \ldots, x_m) = \sum_{x_{m+1} \in G_m}
p_{X_1,\ldots,X_m, X_{m+1}} (x_1, \ldots, x_m, x_{m+1})
\end{equation}
as a direct consequence of total probability.

Among all possible stochastic processes on $G_M$, we will consider
{\em homogeneous Markov chains}. They are fully characterized
by the {\em initial probability}
\begin{equation}
\label{initprob}
p(x) \stackrel{\text{def}}{=} p_{X_1} (x) = \mathbb{P}(X_1 = x)
\end{equation}
and by the transition probability
\begin{equation}
\label{transprob}
P(x,y) \stackrel{\text{def}}{=} \mathbb{P}(X_{m+1} = y|X_m = x)
\end{equation}
that does not depend on the specific value of $m$ (hence the adjective homogeneous).
Note that it is convenient to consider the initial probability as a row vector
with $2^{M(M+1)/2}$ entries with the property that
\begin{equation}
\label{initprobprop}
\sum_{x \in G_M} p(x) = 1,
\end{equation}
and the transition probability as a
$2^{M(M+1)/2} \times 2^{M(M+1)/2}$ matrix, also called {\em stochastic matrix}
with the property that
\begin{equation}
\label{transprobprop}
\sum_{y \in G_M} P(x,y) = 1.
\end{equation}
For a homogeneous Markov chain, the finite dimensional distributions are given by
\begin{equation}
\label{fidimarkov}
p_{X_1,\ldots,X_m} (x_1,\ldots,x_m) = p(x_1) P(x_1,x_2) \cdots P(x_{m-1},x_m).
\end{equation}
It is a well known fact that the finite dimensional distributions
in equation (\ref{fidimarkov}) do satisfy Kolmogorov's conditions
(\ref{sym}) and (\ref{cond2}). Kolmogorov's extension theorem then
implies the existence of Markov chains \cite{billingsley}.
Marginalization of equation (\ref{fidimarkov}) leads to a formula
for $p_{X_m} (x_m) = \mathbb{P} (X_m = x_m)$, this is given by
\begin{equation}
\label{markovgeneral} p_{X_m} (x_m) = \sum_{x_1 \in G_M} p(x_1)
P^{m-1} (x_1,x_m)
\end{equation}
where $P^{m-1}(x_1,x_m)$ is the entry $(x_1,x_m)$ of the
$(m-1)$-th power of the stochastic matrix. Note that, from
equation (\ref{markovgeneral}) and homogeneity one can prove the
{\em Markov semi-group property}
\begin{equation}
\label{semigroup}
P^{m+r}(x,y) = \sum_{z \in G_M} P^m (x,z) P^r (x,y).
\end{equation}

Starting from the basic Markov process with the set of graphs as
space state, we can also consider other auxiliary processes. Just
to mention few among them, we recall:
\begin{itemize}
\item the process counting the number of edges (i.e., the sum of
the adjacency matrix $\mathbb{A}$);

\item the process recording the degree of the graph (i.e., the
marginal total of the adjacency matrix $\mathbb{A}$);

\item the process which measures the cardinality of the strongly
connected components of the graph.
\end{itemize}

Notice that the function of a Markov chain is not a Markov chain
in general, and therefore the study of such processes is not
trivial.

Under a more combinatorial approach, one can consider also the
process recording the permanent of the adjacency matrix
$\mathbb{A}$. We recall that the permanent of the matrix
$\mathbb{A}$ is given by
\begin{equation}
\mathrm{perm}(\mathbb{A})=\sum_{\sigma\in S_M}\prod_{i=1}^M
{\mathbb A}_{i,\sigma(i)}
\end{equation}
where $S_M$ is the symmetric group on the set $\{1, \ldots, M\}$.
The permanent differs from the best known determinant only in the
signs of the permutations. In fact,
\begin{equation}
\mathrm{det}(\mathbb{A})=\sum_{\sigma\in S_M} (-1)^{|\sigma|}
\prod_{i=1}^M {\mathbb A}_{i,\sigma(i)}
\end{equation}
where $|\sigma|$ is the parity of the permutation $\sigma$. Notice
that the permanent is in general harder to compute than the
determinant, as Gaussian elimination cannot be used. However, the
permanent is more appropriate to study the structure of the
graphs. It is known, see for instance \cite{biggs93}, that the
permanent of the adjacency matrix counts the number of the
bijective functions $\phi: V_M \longrightarrow V_M$. The bijective
functions $\phi$ are known in this context as perfect matchings,
i.e., the rearrangements of the vertices consistent with the edges
of the graph. The relations between permanent and perfect
matchings are especially studied in the case of bipartite graphs,
see \cite{schrijver03} for a review of some classical results.

Moreover, we can approach the problem also from the point of view
of {\em symbolic computation}, and we introduce the {\em permanent
polynomial}, defined for each adjacency matrix as follows. Let
${\mathbb Y}$ be an $M \times M$ matrix of variables ${\mathbb Y}
= (y_{i,j})_{i,j=1}^m$. The permanent polynomial is the polynomial
\begin{equation} \label{pperm}
\mathrm{pperm}({\mathbb A}) = \mathrm{perm}({\mathbb Y} \odot
{\mathbb A}) \ ,
\end{equation}
where $\odot$ denotes the element-wise product. For example, the
polynomial determinant of the adjacency matrix
$$
\mathbb{A} = \left( \begin{array}{ccc}
0 & 1 & 1\\
1 & 0 & 1 \\
1 & 1 & 0
\end{array} \right)
$$
introduced above is
$$
\mathrm{pperm}({\mathbb A}) = \mathrm{det} \left(
\begin{array}{ccc}
0 & y_{1,2} & y_{1,3} \\
y_{2,1} & 0 & y_{2,3} \\
y_{3,1} & y_{3,2} & 0
\end{array} \right) = y_{1,2}y_{2,3}y_{3,1} +
y_{1,3}y_{3,2}y_{2,1} \, .
$$
The permanent polynomial in Equation \eqref{pperm} is a
homogeneous polynomial with degree $M$ and it has as many terms as
the permanent of ${\mathbb A}$, all monomials are pure (i.e., with
unitary coefficient) and each transition of the Markov chain from
the adjacency matrix ${\mathbb A}_1$ to the matrix ${\mathbb A}_2$
induces a polynomial $\mathrm{pperm}({\mathbb
A}_2)-\mathrm{pperm}({\mathbb A}_1)$.

Finally, is is also interesting to consider conditional graphs.
With this term we refer to processes on a subset of the whole
family of graphs ${G}_M$. For instance we may require to move only
between graphs with a fixed degree, i.e., between adjacency
matrices with fixed row (and column) totals. In such a case, also
the construction of a connected Markov chain in discrete time is
an open problem, recently approached through algebraic and
combinatorial techniques based on the notion of Markov basis, see
\cite{drton|sturmfels|sullivant09,rapallo06, rapallo|yoshida10}.
This research topic, named {\em Algebraic Statistics} for
contingency tables, seems to be promising when applied to
adjacency matrices of graphs.

\subsection{Ingredient 2: a semi-Markov counting process}

Let $J_1, \ldots, J_n, \ldots$ be a sequence of positive independent and identically distributed (i.i.d.)
random variables interpreted as sojourn times between events in a point process. They are a
{\em renewal process}. Let
\begin{equation}
\label{epochs}
T_n = \sum_{i=1}^n J_i
\end{equation}
be the {\em epoch} (instant) of the $n$-th event. Then, the process $N(t)$ counting the events
occurred up to time $t$ is defined by
\begin{equation}
\label{counting}
N(t) = \max \{n:\, T_n \leq t \}.
\end{equation}

A well-known (and well-studied) counting process is the Poisson process. If $J \sim \exp(\lambda)$, one can prove that
\begin{equation}
\label{poisson}
P(n,t) \stackrel{\text{def}}{=} \mathbb{P}(N(t) = n) = \exp(-\lambda t) \frac{(\lambda t)^n}{n!}.
\end{equation}
The proof leading to the exponential distribution of sojourn times to the Poisson distribution of the counting
process is rather straightforward. First of all one notices that the event $\{ N(t) < n+1 \}$ is given by
the union of two disjoint events
\begin{equation}
\label{disjoint}
\{ N(t) < n+1 \} = \{N(t) < n\} \cup \{N(t) = n\},
\end{equation}
therefore, one has
\begin{equation}
\label{prob1}
\mathbb{P}(N(t) =n) = \mathbb{P}(N(t) < n+1) - \mathbb{P}(N(t) < n);
\end{equation}
but, by definition, the event $\{ N(t) < n\}$ coincides with the event $\{T_n > t\}$. Therefore, from
equation (\ref{prob1}), one derives that
\begin{equation}
\label{prob2}
\mathbb{P}(N(t) = n) = \mathbb{P}(T_n \leq t) - \mathbb{P}(T_{n+1} \leq t).
\end{equation}
The thesis follows from equation (\ref{epochs}). The cumulative
distribution function of $T_n$ is the $n$-fold convolution of an
exponential distribution, leading to the {\em Erlang (or Gamma)
distribution}
\begin{equation}
\label{erlang}
\mathbb{P}(T_n \leq t) = 1 - \sum_{k=0}^{n-1} \exp(-\lambda t) \frac{(\lambda t)^k}{k!},
\end{equation}
and, by virtue of equation (\ref{prob2}), the difference $\mathbb{P}(T_n \leq t) - \mathbb{P}(T_{n+1} \leq t)$
gives the Poisson distribution of equation (\ref{poisson}). Incidentally, it can be proved that $N(t)$ has stationary
and independent increments.

One can also start from the Poisson
process and then show that the sojourn times are i.i.d. random variables. The Poisson process can be defined as
a non-negative-integer-valued stochastic process $N(t)$ with $N(0) = 0$ and with stationary and independent increments
(i.e. a L\'evy process; well, it must also be stochastically continuous,
that is it must be true that for all $a >0$, and for all $s \geq 0$
$\lim_{t \to s} \mathbb{P}(|N(t) - N(s)| > a) = 0$) such that its increment $N(t)-N(s)$ with $(0 \leq s < t)$
has the following distribution for $n \geq 0$
\begin{equation}
\label{poissonprocess}
\mathbb{P}(N(t) - N(s)=n) = \exp(-\lambda (t-s)) \frac{[\lambda(t-s)]^n}{n!}.
\end{equation}
Based on the definition of the process, it is possible to build any of its finite dimensional distributions
using the increment distribution. For instance $\mathbb{P}(N(t_1) = n_1, N(t_2) = n_2)$ with $t_2 > t_1$ is
given by
\begin{eqnarray}
\label{twopointfidi}
\mathbb{P}(N(t_1)=n_1, N(t_2) = n_2) & = & \nonumber \\
& = & \mathbb{P}(N(t_1) = n_1)\mathbb{P}(N(t_2) - N(t_1) = n_2 - n_1) \nonumber \\
& = & \exp(-\lambda t_1) \frac{(\lambda t_1)^{n_1}}{n_1!} \exp(-\lambda (t_2-t_1))
\frac{[\lambda(t_2-t_1)]^{n_2-n_1}}{(n_2-n_1)!}.
\end{eqnarray}
Every L\'evy process, including the Poisson process is Markovian
and has the so-called {\em strong Markov property} roughly meaning
that the Markov property is true not only for deterministic times,
but also for random stopping times. Using this property, it is
possible to prove that the sojourn times are independent and
identically distributed. For $N(0) = 0$, let $T_n = \inf\{t: N(t)
= n \}$ be the $n$-th epoch of the Poisson process (the time at
which the $n$-th jump takes place) and let $J_k = T_k - T_{k-1}$
be the $k$-th sojourn time ($T_0 = 0$). For what concerns the
identical distribution of sojourn times, one has that
\begin{equation}
\label{sojsurv}
\mathbb{P}(T_1 > t) = \mathbb{P}(J_1 > t) = \mathbb{P}(N(t) = 0) = \exp(-\lambda t),
\end{equation}
and for a generic sojourn time $T_k - T_{k-1}$, one finds
\begin{eqnarray}
\label{sojdistr}
T_k - T_{k+1} & = & \inf \{t-T_{k-1}: N(t) = k\} \nonumber \\
& =& \inf \{t-T_{k-1}: N(t) - N(T_{k-1}) = 1\} \nonumber \\
& \stackrel{\text{d}}{=} & \inf \{t-T_{k-1}: N(t - T_{k-1}) = 1, N(T_{k-1}) = 0 \} \nonumber \\
& = & \inf \{t: N(t)=1, N(0) = 0 \},
\end{eqnarray}
where $\stackrel{\text{d}}{=}$ denotes equality in distribution and the equalities are direct consequences
of the properties defining the Poisson process. The chain of equalities means that every sojourn time
has the same distribution of $J_1$ whose survival function is given in equation (\ref{sojsurv}).
As mentioned above, the independence of sojourn times is due to the strong Markov property.
As a final remark, in this digression on the Poisson process, it is important to notice that
one has that its {\em renewal function} $H(t) \stackrel{\text{def}}{=} \mathbb{E}(N(t))$ is
given by
\begin{equation}
\label{poirenewalfunction}
H(t) = \lambda t
\end{equation}
i.e. the renewal function of the Poisson process linearly grows with time, whereas its renewal
density $h(t)$ defined as
\begin{equation}
\label{defrenewaldensity}
h(t) \stackrel{\text{def}}{=} \frac{dH(t)}{dt}
\end{equation}
is constant:
\begin{equation}
\label{poirenewaldensity}
h(t) = \lambda.
\end{equation}

Here, for the sake of simplicity, we shall only consider renewal processes and the related counting processes
(see equations (\ref{epochs}) and (\ref{counting})). When sojourn times are non-exponentially distributed,
the corresponding counting process $N(t)$ is no longer L\'evy and Markovian, but it belongs
to the class of {\em semi-Markov} processes further characterized in the next section
\cite{cinlar75,flomenbom05,flomenbom07,janssen07}. If $\psi_J (t)$ denotes the
probability density function of sojourn times and $\Psi_J(t) \stackrel{\text{def}}{=} \mathbb{P}(J > t)$ is the
corresponding survival function, it is possible to prove the {\em first renewal equation}
\begin{equation}
\label{firstrenewalequation}
H(t) = 1 - \Psi_J (t) + \int_{0}^t H(t-u) \psi_J (u) \, du,
\end{equation}
as well as the {\em second renewal equation}
\begin{equation}
\label{secondrenewalequation}
h(t) = \psi_J (t) + \int_{0}^t h(t-u) \psi_J (u) \, du.
\end{equation}
The second renewal equation is an immediate consequence of the first one, based on the definition of the renewal density
$h(t)$ and on the fact that $\psi_J(t) = - d\Psi_J (t) / dt$. The first renewal equation can be obtained
from equation (\ref{prob2}) which is valid in general and not only for exponential waiting times.
One has the following chain of equalities
\begin{eqnarray}
\label{renewalproof1}
H(t) & =  & \mathbb{E}(N(t)) = \sum_{n=0}^\infty n \mathbb{P}(N(t) = n) = \sum_{n=0}^\infty
n (\mathbb{P}(T_n \leq t) - \mathbb{P}(T_{n+1} \leq t)) \nonumber \\
& = & \sum_{n=1}^\infty \mathbb{P}(T_n \leq t) = \sum_{n=1}^\infty F_n (t),
\end{eqnarray}
where $F_n (t)$ is the cumulative distribution function of the random variable $T_n$, a sum of
i.i.d. positive random variables. Let $f_n(t)$ represent the corresponding density function.
By taking the Laplace transform of equation (\ref{renewalproof1}) and using the fact that
\begin{equation}
\label{densitysum}
\widetilde{f}_n (s) = [\widetilde{\psi}_J (s)]^n,
\end{equation}
one eventually gets
\begin{equation}
\label{renewalproof2}
\widetilde{H}(s) = \sum_{n=1}^\infty \widetilde{F}_n (s) = \frac{1}{s} \sum_{n=1}^\infty \widetilde{f}_n (s) =
\frac{1}{s} \sum_{n=1}^\infty [\widetilde{\psi}_J (s)]^n = \frac{\widetilde{\psi}_J (s)}{s}
\sum_{m=0}^\infty [\widetilde{\psi}_J (s)]^m = \frac{\widetilde{\psi}_J(s)}{s} \frac{1}{1-\widetilde{\psi}_J(s)},
\end{equation}
or (as $|\widetilde{\psi}_J (s)| < 1$ for $s \neq 0$)
\begin{equation}
\label{renewalproof3}
(1-\psi_J(s))\widetilde{H}(s) = \frac{\widetilde{\psi}_J (s)}{s};
\end{equation}
the inversion of equation (\ref{renewalproof3}) yields the first renewal equation (\ref{firstrenewalequation}).

If the sojourn times have a finite first moment (i.e. $\mu_J = \mathbb{E}(J) < \infty$), one has a
{\em strong law of large numbers} for renewal processes
\begin{equation}
\label{stronglaw}
\lim_{t \to \infty} \frac{N(t)}{t} = \frac{1}{\mu_J}, \; \mathrm{a.s.}
\end{equation}
and as a consequence of this result, one can prove the so-called {\em elementary renewal theorem}
\begin{equation}
\label{elementaryrenewaltheorem}
\lim_{t \to \infty} \frac{H(t)}{t} = \frac{1}{\mu_J}.
\end{equation}
The intuitive meaning of these theorems is as follows: if a renewal process is observed a long time after its inception, it is
impossible to distinguish it from a Poisson process. As mentioned in section \ref{intro}, the elementary renewal theorem can explain the ubiquity
of the Poisson process. After a trasient period, most renewal processes behave as the Poisson process.
However, there is a class of renewal processes for which the condition $\mathbb{E}(J) < \infty$
is not fulfilled. These processes never behave as the Poisson process. A prototypical example is
given by the renewal process of Mittag-Leffler type introduced by one of us together with F.
Mainardi and R. Gorenflo back in 2004 \cite{scalas04,mainardi04}. A detailed description of this
process will be given in one of the examples below.

\subsection{Putting the ingredients together}

Let $X_1, \ldots, X_n$ represent a (finite) Markov chain on the state space $G_M$, we now introduce
the process $Y(t)$ defined as follows
\begin{equation}
\label{subordination}
Y(t) \stackrel{\text{def}}{=} X_{N(t)},
\end{equation}
that is the Markov chain $X_n$ is {\em subordinated} to a counting process $N(t)$ coming from a renewal process
as discussed in the previous subsection, with $X_n$ independent of
$N(t)$. In other words, $Y(t)$ coincides with the Markov chain, but
the number of transitions up to time $t$ is a random variable ruled by the probability law of $N(t)$ and
the sojourn times in each state follow the law characterized by the probability density function
$\psi_J (t)$, or, more generally, by the survival function $\Psi_J (t)$.

As already discussed, such a process belongs to the class of semi-Markov processes
\cite{cinlar75,flomenbom05,flomenbom07,janssen07,germano09}, i.e.\ for any $A \subset
G_M$ and $t > 0$ we do have
\begin{equation}
\label{semi-markov}
\mathbb{P}(X_n \in A, J_n \leq t \,|\, X_0, \ldots, X_{n-1}, J_1, \ldots,
J_{n-1}) \\
= \mathbb{P}(X_n \in A, J_n \leq t \,|\, X_{n-1})
\end{equation}
and, if the state $X_{n-1} = x$ is fixed at time $t_{n-1}$, the
probability on the right-hand side will be independent of $n$.
Indeed, by definition, given the independence between the Markov
chain and the counting process, one can write
\begin{eqnarray}
\label{semi-markov-bis}
\mathbb{P}(X_n \in A, J_n \leq t \,|\, X_0, \ldots, X_{n-1}, J_1, \ldots,
J_{n-1}) & = & \mathbb{P}(X_n \in A \,|\, X_{n-1} = x) \mathbb{P}(J_n \leq t) \nonumber \\
& = & P(x,A) (1 - \Psi_J (t)),
\end{eqnarray}
where
\begin{equation}
\label{transitiontoA}
P(x,A) = \sum_{y \in A} P(x,y).
\end{equation}
Equation (\ref{semi-markov-bis}) is a particular case of (\ref{semi-markov}).

It is possible to introduce a slight
complication and still preserve the semi-Markov property. One can imagine that the sojourn time in each state is
a function of the state itself. In this case $\mathbb{P}(J_n \leq t)$ is no longer independent of the
state of the random variable $X_{n-1}$ and equation (\ref{semi-markov-bis}) is replaced by
\begin{eqnarray}
\label{semi-markov-ter}
\mathbb{P}(X_n \in A, J_n \leq t \,|\, X_0, \ldots, X_{n-1}, J_1, \ldots,
J_{n-1}) & = & \mathbb{P}(X_n \in A \,|\, X_{n-1} = x) \mathbb{P}(J_n \leq t|X_{n-1}=x) \nonumber \\
& = & P(x,A) (1 - \Psi^x_J (t)),
\end{eqnarray}
where $\Psi^x_J (t)$ denotes the state-dependent survival function. However, in this case, the random
variable $T_n$ is still the sum of independent random variables, but they are
no-longer identically distributed, and the analysis of the previous section
has to be modified in order to take this fact into proper account.

\section{Examples}

In order to show the behavior of the stochastic processes
described in the previous sections we have simulated the
distribution of two stopping times in two different situations.
The simulations have been written in {\tt R}, see
\cite{rproject10} and the source files are available as
supplementary online material. Notice that some specific packages
for the analysis of graph structures are available, see for
instance \cite{rgraph10}. However, we have used only the {\tt
R}-base commands, as our examples can be analyzed easily without
any additional package.

The examples in this section are designed to introduce the reader
to the simulation algorithms in a framework as simple as possible.
An extended example about a model of interbank market will be
discussed in the next section.

In our examples we use the Mittag-Leffler distribution for the
sojourn times. We recall that the Mittag-Leffler distribution has
survival function given by
\begin{equation}
\label{mlsurvival}
\Psi_J (t)= \mathbb{P}(J > t) = E_\beta (-t^{\beta}),
\end{equation}
where $E_\beta (z)$ is the one-parameter Mittag-Leffler function defined as
\begin{equation}
\label{mlfunction}
E_\beta (z) = \sum_{n=0}^{\infty} \frac{z^n}{\Gamma(n \beta + 1)},
\end{equation}
for $0 < \beta \leq 1$. There are two strong reasons for this
choice. The first one is that many analytical results are
available on the Mittag-Leffler renewal process a.k.a. fractional
Poisson process
\cite{scalas04,mainardi04,fulger08,beghin09,meerschaert10}. The
second reason is that the Mittag-Leffler distribution is the
repeated-thinning limit of heavy-tailed sojourn-time distributions
with algebraically decaying tails with exponent $0 < \beta < 1$
\cite{mainardi04}. For $\beta =1$, the exponential distribution is
recovered from \eqref{mlsurvival}.

\subsection{First example}

In this example we consider graphs without self-loops. Let us
consider a fixed number $M$ of vertices and define a process as
follows:
\begin{itemize}
\item At the time $0$, there are no edges in the graph;

\item At each time, we choose an edge $e$ with uniform
distribution on the $2^{\frac {M(M-1)} {2}}$ edges. If $e$ belongs
to the graph we remove it; if $e$ does not belong to the graph we
add it;

\item The stopping time is defined as the first time for which a
triangle appears in the graph.
\end{itemize}

To simulate the distribution of the stopping times we have used
$10,000$ replications. As the Mittag-Leffler distribution is
heavy-tailed, the density plot and the empirical distribution
function plot are not informative. Thus, we have reported the
box-plot, to highlight the influence of the outliers.

\begin{figure}
\epsfig{file=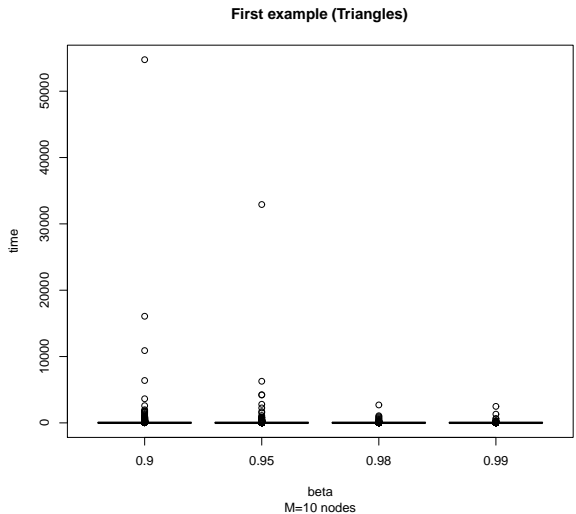, width=80mm} \caption{Box-plot of the
distribution of the stopping times with varying $\beta$ for
Example A.} \label{figureex1_1}
\end{figure}

With a first experiment, we have studied the influence of the
$\beta$ parameter. In graph with $M=10$ nodes, we have considered
the sojourn times with a Mittag-Leffler distribution with
different $\beta$ parameter, namely $\beta=0.90. 0.95, 0.98,
0.99$. The box-plot are displayed in Figure \ref{figureex1_1}, and
some numerical indices are in Table \ref{tableex1}.

\begin{table}
\begin{center}
\begin{footnotesize}
\begin{verbatim}
beta: 0.9
     Min.   1st Qu.    Median      Mean   3rd Qu.      Max.
3.628e-01 8.010e+00 1.275e+01 3.146e+01 2.026e+01 5.475e+04
------------------------------------------------------------
beta: 0.95
     Min.   1st Qu.    Median      Mean   3rd Qu.      Max.
3.249e-01 7.545e+00 1.144e+01 2.086e+01 1.689e+01 3.292e+04
------------------------------------------------------------
beta: 0.98
     Min.   1st Qu.    Median      Mean   3rd Qu.      Max.
   0.2607    7.2960   10.8600   12.9500   15.0500 2704.0000
------------------------------------------------------------
beta: 0.99
     Min.   1st Qu.    Median      Mean   3rd Qu.      Max.
   0.5373    7.2190   10.6300   12.3400   14.6700 2487.0000
\end{verbatim}
\end{footnotesize}
\caption{Summary statistics for Example A with varying $\beta$.} \label{tableex1}
\end{center}
\end{table}

Our results show that:
\begin{itemize}
\item the outliers are highly influenced from the value of
$\beta$. This holds, with a less strong evidence, also for the
quartiles $Q1$ and $Q3$;

\item the median is near constant, while the mean is affected by
the presence of outliers.
\end{itemize}

\begin{figure}
\begin{tabular}{cc}
\epsfig{file=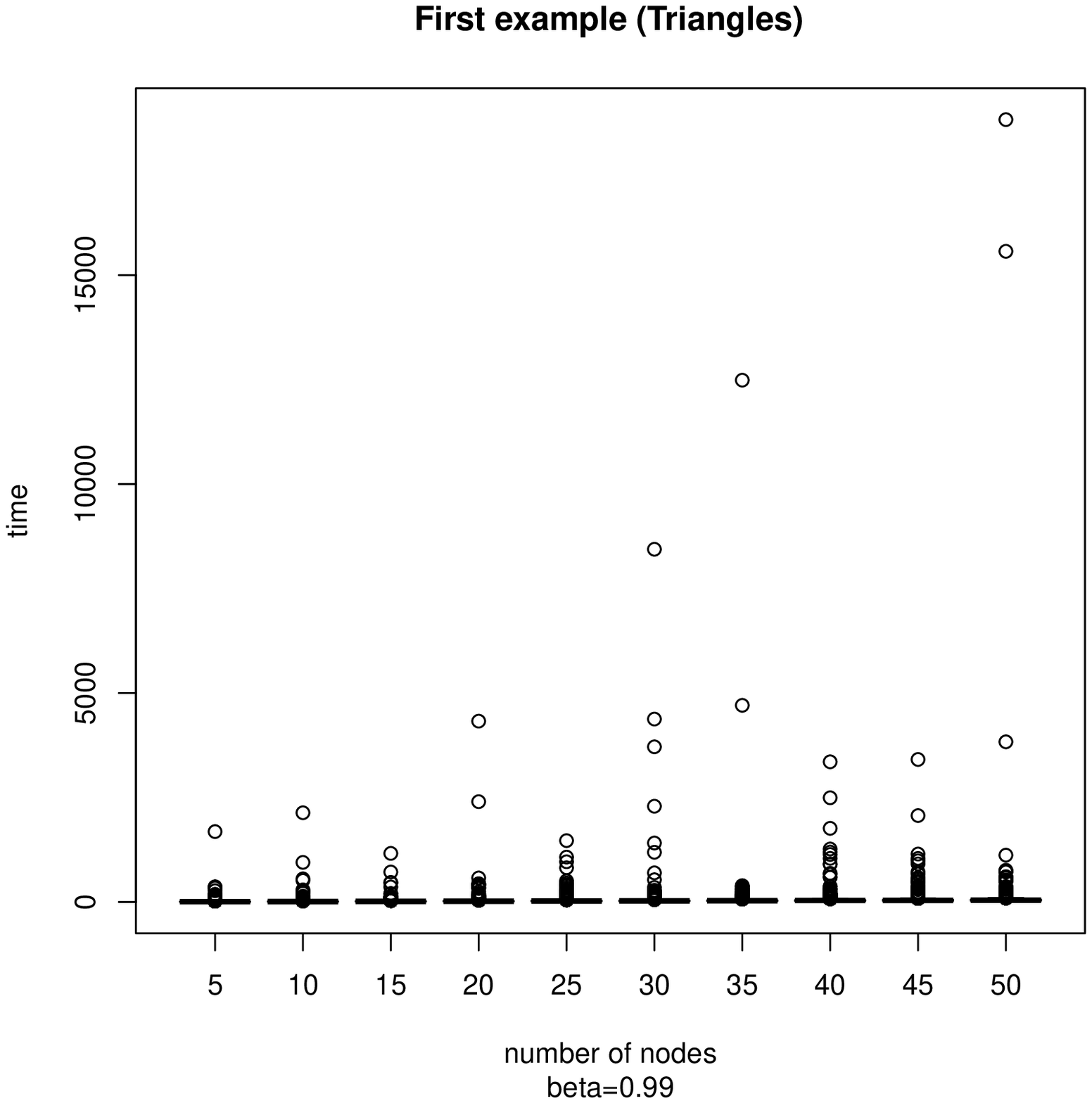, width=80mm} & \epsfig{file=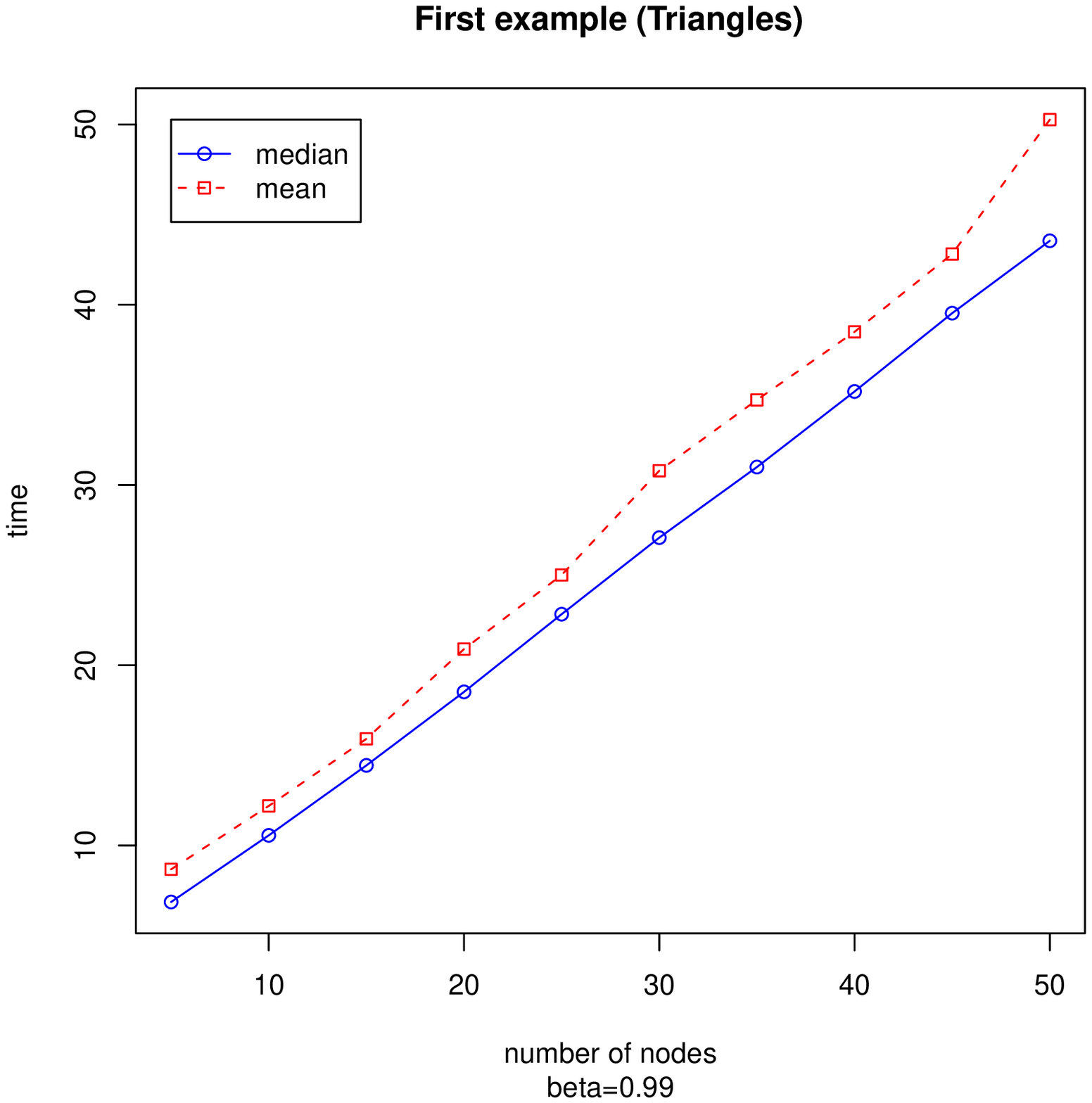,
width=80mm}
\end{tabular}
\caption{Box-plot (left), mean and median (right) of the
distribution of the stopping times with varying $M$ for Example
A.} \label{figurex1_2}
\end{figure}

With a second experiment, we have considered a fixed parameter
$\beta=0.99$, but a variable number of vertices $M$ ranging from
$5$ to $50$ by $5$. In Figure \ref{figurex1_2} we present the
box-plots of the stopping time distribution and the trends of the
mean and the median.

From this graphs we can notice that:
\begin{itemize}
\item the presence of outliers is more relevant in the graph with
a large number of nodes;

\item the mean and the median are roughly linear, but the trend of
the median looks more stable.
\end{itemize}

\subsection{Second example}

Let us consider a population with individuals $\{1, \ldots , M \}$
and suppose that the person labelled 1 has to share some
information. At a first random time, he chooses another individual
with random uniform probability and shares the information with
him. At a second random time, one person who has the information
chooses an individual among the other $(M-1)$ and shares again the
information. Note that each individual shares the information with
another one, no matters if he has already or not the information.
In terms of graphs, we define a process as follows:

\begin{itemize}
\item At the time $0$, there are no edges in the graph;

\item At each time, we choose a vertex $m$ connected with 1 and we
add to the graph an edge among $(m,1), (m,2), \ldots,
(m,m-1),(m,m+1), \ldots, (m,M-1),(m,M)$ with random uniform
distribution. If the chosen edge is already in the graph we do
nothing;

\item The stopping time is defined as the first time for which the
whole graph is connected.
\end{itemize}

\begin{figure}
\epsfig{file=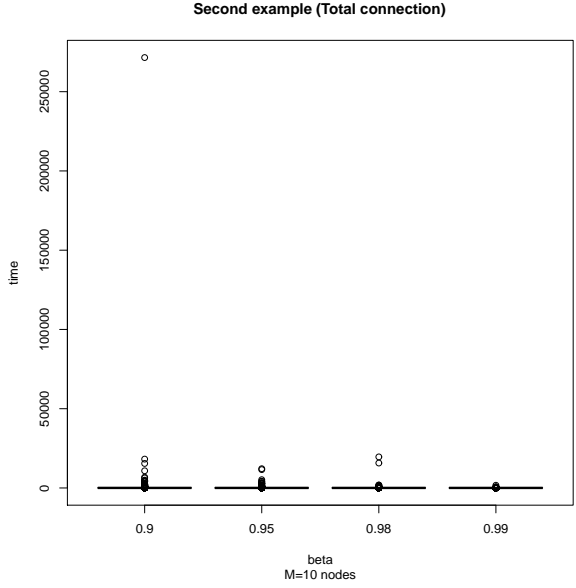, width=80mm} \caption{Box-plot of the
distribution of the stopping times with varying $\beta$ for
Example B.} \label{figureex2_1}
\end{figure}

The experimental settings for this example are the same as for
Example A. With a fixed number of vertices $M=10$ and varying
$\beta$ as above, we obtain the box-plots in Figure
\ref{figureex2_1}, and the numerical summary in Table
\ref{tableex2}.

\begin{table}
\begin{center}
\begin{footnotesize}

\begin{verbatim}
beta: 0.9
     Min.   1st Qu.    Median      Mean   3rd Qu.      Max.
3.786e+00 2.154e+01 3.207e+01 8.811e+01 4.938e+01 2.715e+05
------------------------------------------------------------
beta: 0.95
     Min.   1st Qu.    Median      Mean   3rd Qu.      Max.
    3.393    19.410    27.050    41.550    38.260 12140.000
------------------------------------------------------------
beta: 0.98
     Min.   1st Qu.    Median      Mean   3rd Qu.      Max.
    3.565    18.230    24.980    33.530    33.970 19600.000
------------------------------------------------------------
beta: 0.99
     Min.   1st Qu.    Median      Mean   3rd Qu.      Max.
    4.738    17.690    23.940    27.160    32.310  1701.000
\end{verbatim}

\end{footnotesize}
\end{center}
\caption{Summary statistics for Example B with varying $\beta$.} \label{tableex2}

\end{table}

From this results we can see that the outliers are highly
influenced from the value of $\beta$, while the variation of the quantiles $Q1$ and
$Q3$ is much lower. Also in this example, the mean is
affected by the presence of outliers.

\begin{figure}
\begin{tabular}{cc}
\epsfig{file=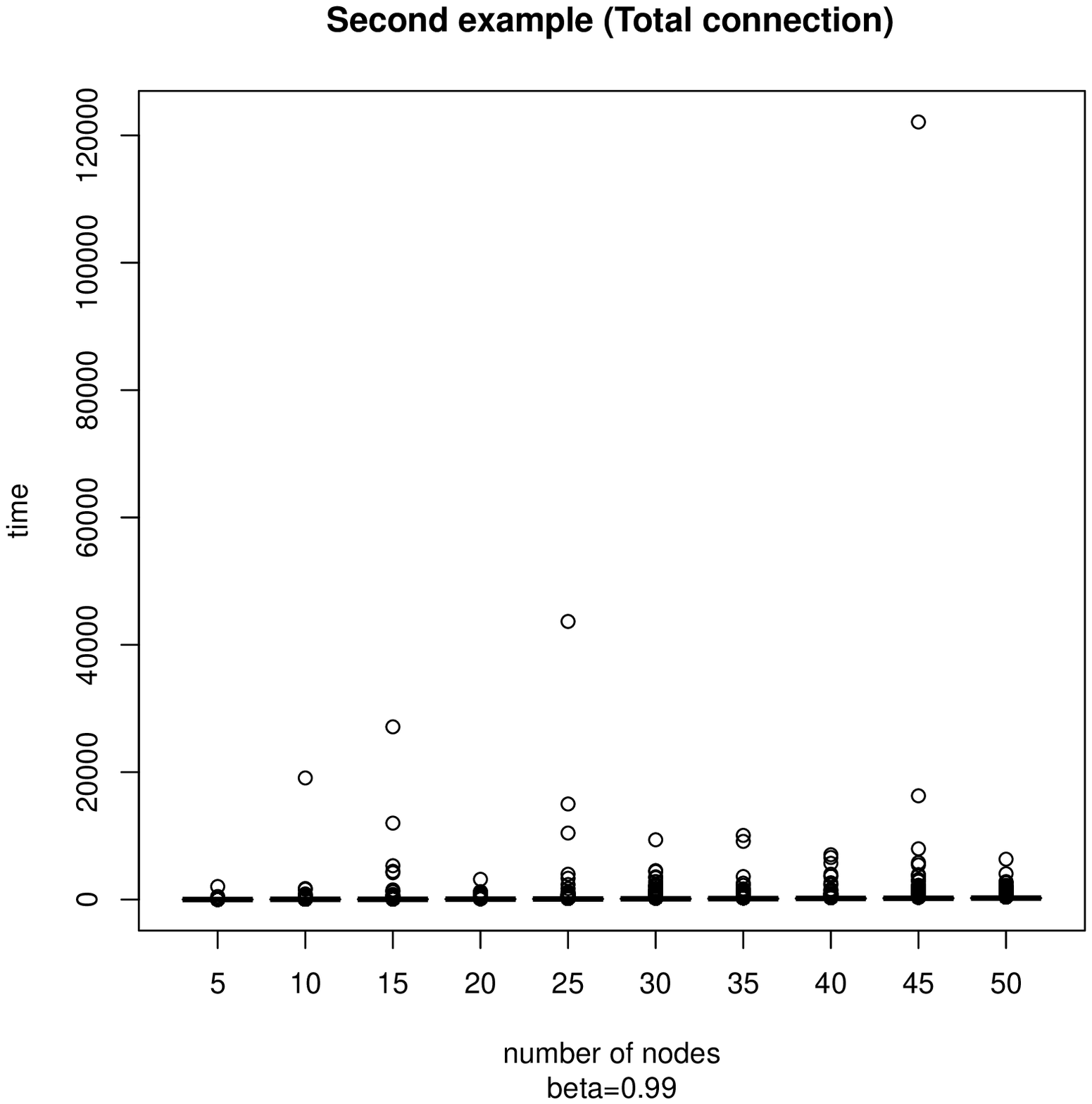, width=80mm} & \epsfig{file=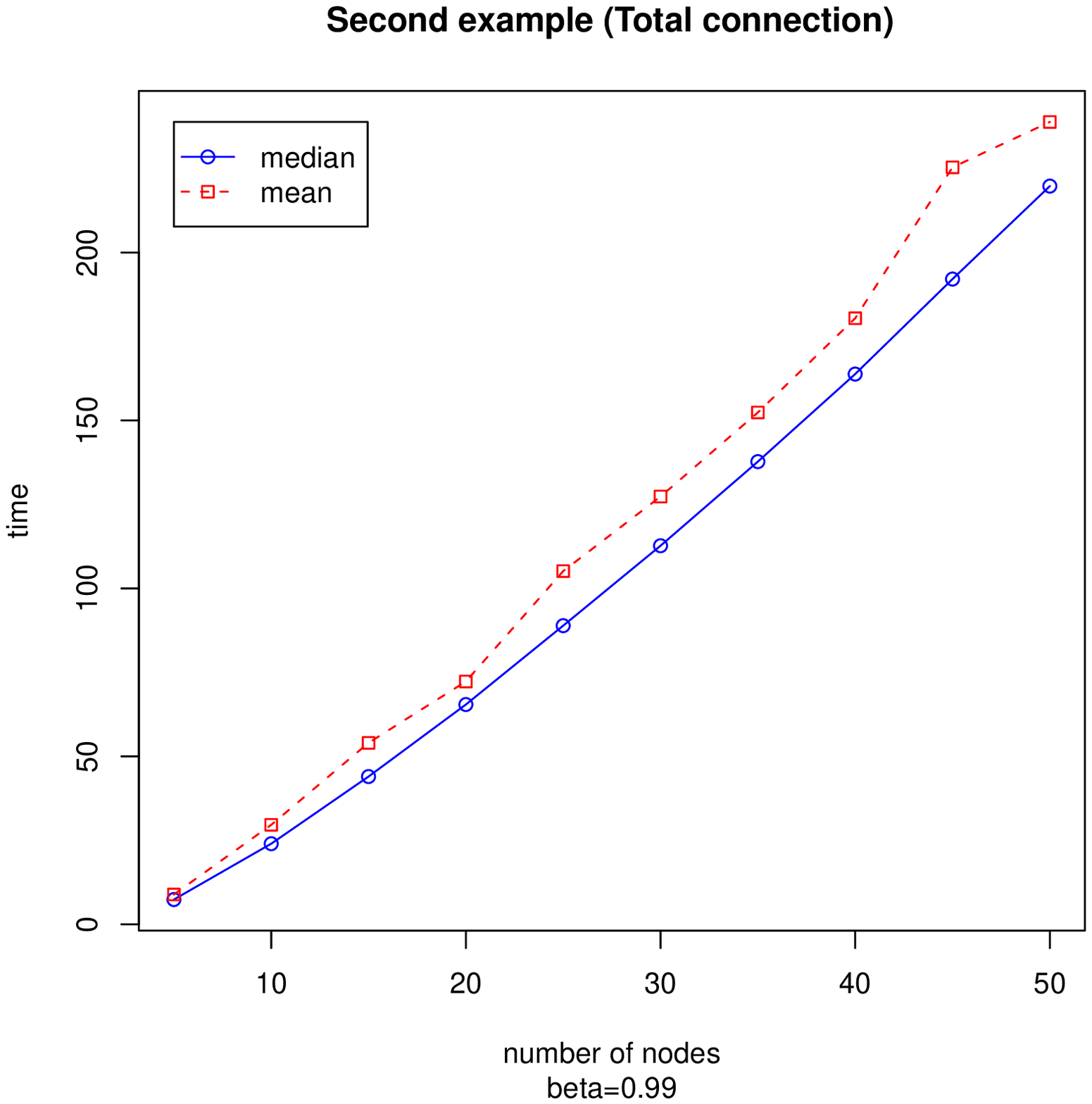,
width=80mm}
\end{tabular}
\caption{Box-plot (left), mean and median (right) of the
distribution of the stopping times with varying $M$ for Example
B.} \label{figurex2_2}
\end{figure}

With the second experiment with a variable number of vertices $M$
ranging from $5$ to $50$ by $5$, we obtain the plots displayed in
Figure \ref{figurex2_2}. The conclusions are the same as in the
previous example.

\section{Extended example. An interbank market}

In this section we present a simple model for interbank markets.
It serves the purpose of illustrating the modelling potentialities
of the ideas presented above.

This example deals with an interbank market receiving loan
requests from the corporate sector at random times. For the sake
of readability, in this section we will use the symbol $\Delta
t_k$ instead of $J_k$ for the $k$-th inter-arrival duration and we
will denote the epochs at which loan requests are made with the
symbol $t_k$ instead of $T_k$. In this more realistic example, we
will briefly discuss the difficulties that must be faced when one
desires to go beyond a mere phenomenological description of
reality.

\begin{table}[h]
  \centering
  \begin{tabular}{|l|l|}
  \hline
 \textbf{Assets}  &  \textbf{Liabilities} \\ \hline

   $C^b_{t_n}$: liquidity  & $D^b_{t_n}$: total (households' and firms') deposits \\

   $L^b_{t_n}$: loans to the corporate sector & $B^b_{t_n}$: debt with other banks  \\

    $\mathcal{L}^b_{t_n}$: loans to other banks  & $E^b_{t_n}$: equity (net worth)\\

  \hline
\end{tabular}
\caption{Balance sheet entries of bank $b$ at time $t_n$}
\label{bsbank}
\end{table}

We consider an interbank market characterized by $M$ banks that demand and supply liquidity at a given interest rate $r_B$. Each bank $b$ is described at any time by its balance sheet, as outlined in Table \ref{bsbank}. The market is decentralized and banks exchange liquidity by means of pairwise interactions. Banks lend money also to the corporate sector at the constant rate $r_C > r_B$ and all corporate and interbank loans are to be repayed after $T$ units of time. We stipulate that the loan requests from the corporate sector to the banking system are the events triggering the interbank market and we model these events as a Poisson process of parameter $\lambda$.  In particular, we state that, at exponentially distributed intervals of time $\Delta t_n = t_n - t_{n-1}$, a loan request of constant amount $\ell$ is submitted from the corporate sector to a bank chosen at random with uniform distribution among the $M$ banks. As in the previous examples, in principle, the Poisson process can be replaced by any suitable counting process. Let us denote the chosen bank with the index $i$ and the time at which the loan is requested as $t_n$. If $C^i_{t_{n-1}} < \ell$, the chosen bank is short of liquidity to grant the entire amount of the loan. Given the interest rate spread between $r_C$ and $r_B$, the profit-seeking bank enters the interbank market in order to borrow at the rate $r_B$ the amount $\ell - C^i_{t_{n-1}}$ necessary to grant the entire loan.
In the interbank market, a new bank is then chosen at random with uniform distribution among the remaining $M-1$ banks. Let us denote with $j$ the new chosen bank. If bank $j$ has not enough liquidity to lend the requested amount, i.e., if $C_{t_{n-1}}^j < \ell - C_{t_{n-1}}^i$, then a new bank $h$ is again chosen at random among the remaining $M-2$ ones to provide the residual liquidity, and so on. This process in the interbank market continues till the liquidity amount $\ell - C_{t_{n-1}}^i$ necessary to bank $i$ is collected.

Finally, as soon as the loan $\ell$ is provided to the corporate sector, we stipulate that the deposits as well as the liquidity of any bank $b$, being $b = 1,\ldots,M$, is increased by the amount $\omega_t^b \, \ell$, where $\omega_t^b$ are random numbers constrained by $\sum_{b} \omega_t^b = 1$. The rationale behind this choice is that a loan, when is taken and spent, creates a deposit in the bank account of the agent to whom the payment is made; for instance, when the corporate sector gets a loan to pay wages to workers or to pay investments to capital goods producers, then the deposits at the $M$ banks of the agents receiving the money are increased by a fraction of the borrowed amount $\ell$. We assume that these deposits are randomly distributed among the $M$ banks.

To give a clearer idea on how the balance sheets of banks evolve after an event in the interbank market, let us consider an example where at time $t_n$ the corporate sector requests a loan $\ell$ to the randomly selected bank $i$, which, being short of liquidity (i.e. $C_{t_{n-1}}^i < \ell$), needs to enter into interbank market where it borrows a loan of amount $\ell - C_{t_{n-1}}^i$ from the randomly selected bank $j$. We suppose here $C_{t_{n-1}}^j > \ell - C_{t_{n-1}}^i$, therefore no other lending banks enter the interbank market.
According to the model outlined above, at the end of the interbank market session, the balance sheets of bank $i$ and of bank $j$ change as outlined in Table \ref{bsbank_dynamics}.

\begin{table}[h]
  \centering
  \begin{tabular}{|rcl|rcl|}
  \hline
  $C_{t_n}^i$ & = & $\omega_{t_n}^i \, \ell$ & $C_{t_n}^j$ & = & $C^j_{t_{n-1}} - (\ell-C_{t_{n-1}}^i) + \omega_{t_n}^j \, \ell$ \\
  $L_{t_n}^i$ &=& $L_{t_{n-1}}^i + \ell$ & $L_{t_n}^j$ &=& $L_{t_{n-1}}^j$\\
$\mathcal{L}_{t_n}^i$ &=& $\mathcal{L}_{t_{n-1}}^i$ & $\mathcal{L}_{t_n}^j$ &=& $\mathcal{L}_{t_{n-1}}^j + (\ell-C_{t_{n-1}}^i)$\\
$D_{t_n}^i$ &=& $D_{t_{n-1}}^i + \omega_{t_n}^i \, \ell$ & $D_{t_n}^j$ &=& $D_{t_{n-1}}^j + \omega_{t_n}^j \, \ell$\\
$B_{t_n}^i$ &=& $B_{t_{n-1}}^i + (\ell - C_{t_{n-1}}^i)$ & $B_{t_n}^j$ &=& $B_{t_{n-1}}^j$\\
$E_{t_n}^i$ &=& $E_{t_{n-1}}^i$ & $E_{t_n}^j$ &=& $E_{t_{n-1}}^j$\\
  \hline
\end{tabular}
\caption{Dynamics of balance sheet entries of bank $i$ (lender to the corporate sector and borrower in the interbank market) and bank $j$ (lender in the interbank market) at time $t_n$ when both the corporate loan $\ell$ and the related interbank loan $\ell-C_{t_{n-1}}^i$ are granted.}
\label{bsbank_dynamics}
\end{table}
Once the assets and the debt liabilities entries of any bank are updated following the lending activity to the corporate sector and the interbank market outcomes, the equity is then updated as residual according to the usual accounting equation:
\begin{equation}\label{Eq:equity}
    E^b_{t_n} = C^b_{t_n} + L^b_{t_n} + \mathcal{L}^b_{t_n} - D^b_{t_n} - B^b_{t_n}\,.
\end{equation}
It is worth noting that, as reported in Table \ref{bsbank_dynamics}, the equity of both bank $i$ and $j$ does not change from $t_{n-1}$ to $t_n$. This result is obtained by computing the new equity levels at time $t_n$ using \eqref{Eq:equity} but should not be a surprise given that lending and borrowing clearly change the balance sheet entries of banks but not their net worth at the time the loan is granted or received. Indeed, the net worth of the lending banks is increased by the interest revenues when the corporate loan as well as the interbank loan is repaid together with the interest amounts. In particular, equity of bank $i$ is increased by $r_C \, \ell - r_B \, (\ell-C^i_{t_{n-1}})$, while equity of bank $j$ is increased by $r_B \, (\ell-C^i_{t_{n-1}})$. Table \ref{bsbank_dynamics_repayment} shows how balance sheet entries change at time $t_m = t_n + T$ when the two loans are paid back. It is worth noting again that the equity dynamics is consistent with the dynamics of other balance sheet entries, according to \eqref{Eq:equity}. Finally, as granting a bank loan to the corporate sector increases private deposits at banks, also the opposite holds when a loan is paid back. The repayment of the loan $\ell$ together with interests $r_C \, \ell$ corresponds to a reduction of private deposits, as well as of the related banks' liquidity, of the same amount. As in the previous case, we assume that the reduction $(1+r_C)\, \ell$ is uniformly and randomly distributed among the $M$ banks with weights $\omega_{t_m}^b\,$, where $b = 1,\ldots,M$.

\begin{table}[h]
  \centering
  \begin{tabular}{|rcl|rcl|}
  \hline
  $C_{t_m}^i$ & = & $C_{t_{m-1}}^i + \, (r_C-r_B)\,\ell + (1+r_B)\,C_{t_{n-1}}^i + $ & $C_{t_m}^j$ & = & $C^j_{t_{m-1}} + (1+r_B)\, (\ell-C_{t_{n-1}}^i) \, + $ \\
   &  & $-\, \omega_{t_m}^i \, (1+r_C)\,\ell$ &  &  & $-\, \omega_{t_m}^j \, (1+r_C)\,\ell$\\
  $L_{t_m}^i$ &=& $L_{t_{m-1}}^i - \ell$ & $L_{t_m}^j$ &=& $L_{t_{m-1}}^j$\\
$\mathcal{L}_{t_n}^i$ &=& $\mathcal{L}_{t_{n-1}}^i$ & $\mathcal{L}_{t_m}^j$ &=& $\mathcal{L}_{t_{m-1}}^j - (\ell-C_{t_{n-1}}^i)$\\
$D_{t_m}^i$ &=& $D_{t_{m-1}}^i - \omega_{t_m}^i \, (1+r_C)\, \ell$ & $D_{t_m}^j$ &=& $D_{t_{m-1}}^j - \omega_{t_m}^j \, (1+r_C)\,\ell$\\
$B_{t_m}^i$ &=& $B_{t_{m-1}}^i - (\ell - C_{t_{n-1}}^i)$ & $B_{t_m}^j$ &=& $B_{t_{m-1}}^j$\\
$E_{t_m}^i$ &=& $E_{t_{m-1}}^i + \, (r_C-r_B)\,\ell + r_B\,C_{t_{n-1}}^i $ & $E_{t_m}^j$ &=& $E_{t_{m-1}}^j + r_B \, (\ell-C_{t_{n-1}}^i) $\\
  \hline
\end{tabular}
\caption{Dynamics of balance sheet entries of bank $i$ (lender to the corporate sector and borrower in the interbank market) and bank $j$ (lender in the interbank market) at time $t_m = t_n + T$ when both the corporate loan $\ell$ and the related interbank loan $\ell-C_{t_{n-1}}^i$ are paid back.}
\label{bsbank_dynamics_repayment}
\end{table}

We can then define a $M \times M$ adjacency matrix $\mathbb{A}$
representing the graph associated to the interbank market, where
the nodes of the graph correspond the $M$ banks  and the edges to
the lending and borrowing relationships among banks. At variance with the previous discussion
and examples, here, it is meaningful to consider directed graphs
and therefore the matrix can be asymmetric. In particular, if bank
$j$ is lending money to bank $i$, we set $A_{j,i} = 1$, but we may
have $A_{i,j} = 1$ or $A_{i,j} = 0$, depending if bank $i$ is
lending or not money to bank $j$. The situation where both
$A_{j,i}$ and $A_{i,j}$ are set to 1 is not contradictory but it
means that two loans have been granted in the two opposite
directions, i.e. from bank $i$ to bank $j$ and from bank $j$ to
bank $i$, at different times. The time evolution of the adjacency
matrix depends  on the evolution of events in the interbank
market. In particular, when the first loan from bank $j$ to bank
$i$ is paid back, $A_{j,i}$ is again set to 0, provided that no
new loans have been granted by bank $j$ to bank $i$ in the
meantime, if this happens the value of $A_{j,i}$ remains at 1 till
there are debts of bank $i$ to bank $j$. If this is required by
the particular application, it is even possible to consider {\em
weighted graphs} where the entry $A_{i,j}$ contains the value of
the loan from bank $i$ to bank $j$.

The dynamics in the interbank market can then be represented as a
Markov chain on graphs subordinated to the Poisson process
representing the random events of loan requests to the banking
system by the corporate sector. It is worth noting that the Markov
process and the Poisson process are independent here, however, the
transition probabilities of the Markov process are not fixed
\emph{ex ante} but depends on the endogenous evolution of the
balance sheets of banks. Therefore, here, the Markov process is
{\em not homogeneous}.

\section{Concluding considerations}

We have discussed a model of graph dynamics based on two
ingredients. The first ingredient is a Markov chain on the space
of possible graphs. The second ingredient is a semi-Markov
counting process of renewal type. The model consists in
subordinating the Markov chain to the semi-Markov counting
process. In simple words, this means that the chain transitions
occur at random time instants called epochs.
This model takes into account the fact that social interactions are intrinsically volatile and not permanent.

Note that state dependent subordination (see equation \eqref{semi-markov-ter}) gives rise to very interesting dynamics
from the ergodicity viewpoint \cite{saa}. In order to illustrate this fact, let us consider a simple two-state aperiodic and irreducible Markov chain with the following transition probability matrix:
$$
{P} = \left( \begin{array}{ccc}
0.1 & 0.9\\
0.9 & 0.1
\end{array} \right).
$$
In this case, the invariant measure is uniform and it is given by
$$
p = \left( \frac{1}{2}, \frac{1}{2} \right),
$$
meaning that the probability of finding each state at equilibrium is $1/2$. Now, let us call $A$ the first state and $B$ the second state. Let
the sojourn time in $A$ be exponentially distributed with parameter $\lambda_A$ and the sojourn time in $B$ still exponentially
distributed with parameter $\lambda_B$. If a single realization of this process is considered, the average time of permanence
in state $A$ will be given by $1/\lambda_A$ and the average time of permanence in $B$ will be given by $1/\lambda_B$. Therefore,
if $\lambda_A \neq \lambda_B$, then the ratio of average sojourn times will be different from $1$. In other words, for this simple model,
the fraction of sojourn times is not equal to the fraction of the ensemble measure: a signal of non-ergodicity.

Finally, with reference to the examples discussed above, this kind of modeling can be used for risk evaluation. Given a loss function,
a function that gives the losses when adverse events take place, the risk function is defined as the expected value of the loss function.
With our approach, one can derive the probability of the adverse events as a function of time and use this measure to evaluate the
risk function.

\section{Acknowledgements}

FR wishes to acknowledge an INDAM (Istituto Nazionale di Alta
Matematica, Italy) grant with which this work was partially
funded. Long ago (in 2005 and 2006), ES discussed these topics
with Gilles Daniel, Lev Muchnik and Sorin Solomon.  He now
gratefully acknowledges these discussions. ES wishes to thank Elena
Akhmatskaya for recent discussion on these issues. MR and ES first
discussed these issues during a visit of ES to Reykjavik
University sponsored by the ERASMUS program. ES is also grateful
to Universitat Jaume I for the financial support received from
their Research Promotion Plan 2010 during his scientific visit in
Castell\'on de la Plana where parts this paper were written.

\end{document}